\documentclass[a4paper,reqno]{amsart}

\calclayout

\usepackage{graphicx}%
\usepackage{amsmath,amssymb}%
\usepackage{mathtools}%
\usepackage{eucal}%
\usepackage{xcolor}%
\usepackage{tikz}%
\usetikzlibrary{decorations,snakes}%
\usepackage{tikz-cd}%

\usepackage[%
backend=bibtex,%
style=alphabetic,%
natbib=false,%
isbn=false, 
doi=false, 
url=false, 
eprint=true, 
giveninits=true, 
backref=false, 
]{biblatex} \AtEveryBibitem{%
  \ifentrytype{unpublished}{}{\clearfield{note}}%
  } \renewbibmacro{in:}{} 
\usepackage{doi}%

\usepackage{hyperref}%
\hypersetup{%
  colorlinks,%
  linkcolor={red!80!black},%
  citecolor={blue!80!black},%
  urlcolor={blue!80!black}%
}

\usepackage{thmtools}

\declaretheorem[style=theorem,numbered=no]{theorem}%
\declaretheorem[style=remark,numbered=no]{remark}

\newcommand{\kk}{\mathbf{k}}%

\NewDocumentCommand{\A}{}{A_\infty}%
\RenewDocumentCommand{\AA}{}{\mathcal{A}}%
\NewDocumentCommand{\Hom}{O{\AA}D<>{}mm}{\operatorname{Hom}_{#1}^{#2}\!\left(#3,#4\right)}%
\NewDocumentCommand{\Ext}{O{\AA}D<>{*}mm}{\operatorname{Ext}_{#1}^{#2}\!\left(#3,#4\right)}%
\NewDocumentCommand{\RHom}{O{A}mm}{\operatorname{RHom}_{#1}\!\left(#2,#3\right)}%

\RenewDocumentCommand{\mod}{m}{\operatorname{mod}{#1}}%
\NewDocumentCommand{\Db}{m}{\operatorname{D}^{\mathrm{b}}\!\left(#1\right)}%

\NewDocumentCommand{\len}{m}{\operatorname{len}\!\left(#1\right)}%
\NewDocumentCommand{\class}{m}{#1}%

\NewDocumentCommand{\dg}{m}{{#1}_{\mathrm{dg}}}%

\NewDocumentCommand{\F}{}{\mathcal{F}}%

\newcommand{\myvspace}{}

\DeclareSymbolFont{sfoperators}{OT1}{cmss}{m}{n}%
\DeclareSymbolFontAlphabet{\mathsf}{sfoperators}%
\makeatletter%
\def\operator@font{\mathgroup\symsfoperators}%
\makeatother

\title[On a theorem of B.~Keller]{On a theorem of B.~Keller%
  \\on Yoneda algebras of simple modules}%

\author[G.~Jasso]{Gustavo Jasso}%
\address{%
  Lund University, %
  Centre for Mathematical Sciences, %
  Box 118, %
  22100 Lund, %
  Sweden%
}%
\email{gustavo.jasso@math.lu.se}%
\urladdr{https://gustavo.jasso.info}%



\begin{document}

\maketitle

\begin{abstract}
  A theorem of Keller states that the Yoneda algebra of the simple modules over
  a finite-dimensional algebra is generated in cohomological degrees $0$ and $1$
  as a minimal $\A$-algebra. We provide a proof of an extension of Keller's
  theorem to abelian length categories by reducing the problem to a particular
  class of Nakayama algebras, where the claim can be shown by direct
  computation.
\end{abstract}

\sloppy

We work over an arbitrary field $\kk$. Let $A$ be a finite-dimensional algebra
and $S$ the direct sum of a complete set of representatives of the simple
(right) $A$-modules. The Yoneda algebra $\Ext[A]{S}{S}$, as a graded algebra,
does not determine the algebra $A$ up to Morita equivalence, as the following
example shows~(see~\cite[Sec.~2.1]{Kel02a} or~\cite[Ex.~B.2.2]{Mad02}): Let
$A=\kk[x]/(x^\ell)$, $\ell\geq 3$; then
\begin{equation}
  \label{eq:example}
  \Ext[A]{S}{S}\cong\kk[u,v]/(u^2),\qquad |u|=1,\ |v|=2,
\end{equation}
does not depend on $\ell$. On the other hand, the Yoneda algebra of an arbitrary
finite-dimensional algebra inherits, via Kadeishvili's Homotopy Transfer
Theorem~\cite{Kad82,Mar06}, the structure of a minimal\footnote{Recall that an
  $\A$-algebra is minimal if its underlying complex has vanishing differential.}
$\A$-algebra since the Yoneda algebra is the cohomology of the differential
graded algebra $\RHom[A]{S}{S}$. Endowed with this additional $\A$-structure the
Yoneda algebra does determine the algebra $A$ up to Morita
equivalence~\cite[Sec.~7.8]{Kel01}. The purpose of this short article is to
provide a proof of a minor extension of a theorem of Keller~\cite[Sec.~2.2,
Prop.~1(b)]{Kel02a} that is stated below. A proof of Keller's theorem that
utilises the calculus of Massey products was announced by Minamoto
in~\cite{Min16}; our proof should have some similarities with his.
\myvspace{0.5em}

\begin{theorem}[Keller] Let $\AA$ be a $\kk$-linear abelian length
  category~\cite{KV18} with only finitely many pairwise non-isomorphic simple
  objects, for example the category of finite-dimensional modules over a
  finite-dimensional algebra. Let $S_1,\dots,S_n$ be a complete set of
  representatives of the simple objects in $\AA$ and set ${S\coloneqq
    S_1\oplus\cdots\oplus S_n}$. Then, the Yoneda algebra $\Ext{S}{S}$ is
  generated by its homogeneous components of cohomological degrees $0$ and $1$
  as an $\A$-algebra.
\end{theorem}

\begin{remark}
  Since the abelian category $\AA$ is not assumed to have any non-zero
  projective or injective objects, we interpret the Yoneda algebra
  $\Ext[\AA]{S}{S} $ in terms of Yoneda equivalence classes of exact
  sequences~\cite{Yon60}. The results in~\cite[Sec.~III.3]{Ver96} (see in
  particular paragraph~III.3.3.2 therein) readily imply that we can identify the
  Yoneda algebra with the graded algebra $\bigoplus_{k\geq0}\Hom[]{S}{S[k]}$ of
  endomorphisms of $S$ in the derived category of $\AA$. Thus, we may and we
  will identify minimal $\A$-algebra structures on the Yoneda algebra with those
  of the latter graded algebra.
\end{remark}

\begin{remark}
  The isomorphism~\eqref{eq:example} shows that the conclusion of the theorem
  may fail if the Yoneda algebra $\Ext{S}{S}$ is considered as a graded algebra
  only.
\end{remark}

\begin{proof}[Proof of the theorem]
  We use freely the theory of $\A$-categories~\cite{Kel01,Lef03}, as well as the
  theory of differential graded (=DG) categories and their derived
  categories~\cite{Kel94,Kel06}. We also assume familiarity with the
  Auslander--Reiten (=AR) theory of Nakayama algebras, see for
  example~\cite[Ch.~V]{ASS06}.

  Let $d\geq1$ and $\class{\delta}\in\Ext<d+1>{S}{S}$ a Yoneda class represented
  by an exact sequence
  \begin{equation}
    \label{eq:delta}
    0\to S_a\to M_1\to M_2\to\cdots\to M_{d+1}\to S_b\to0
  \end{equation}
  between some simple objects; notice that we identify
  $\Ext<k>{S}{S}=\bigoplus_{i,j=1}^n\Ext<k>{S_i}{S_j}$. Below we prove that
  $\delta$ is generated by Yoneda classes of degrees $1$ and $d$ under the
  operations $m_k^{\AA}$ with $k\leq\len{M_1}$. The theorem then follows by
  induction on $d\geq1$.

  A straightforward inductive argument shows that the exact sequence
  \eqref{eq:delta} can be extended, by first choosing composition series of
  $N_i\coloneqq\operatorname{img}(M_i\to M_{i+1})$, to a commutative diagram of
  the following form in which all squares are bicartesian, the apparent
  sequences at the bottom are short exact, and all of the objects along the
  bottom of the diagram are simple:
  \begin{center}
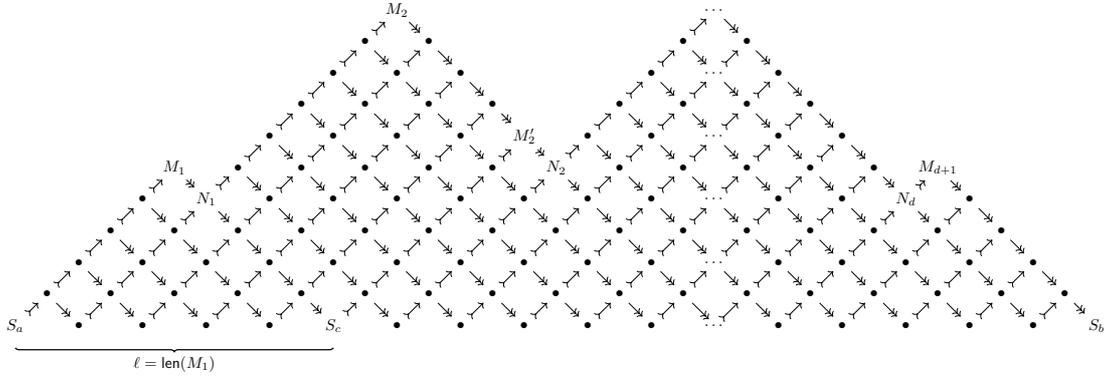
    \begin{resizebox}{\textwidth}{!}{%
    \begin{tikzpicture}[rotate=-45]
      \node (000) at (0,0) {$S_a$};%
      \node (001) at (0,1) {$\bullet$};%
      \node (002) at (0,2) {$\bullet$};%
      \node (003) at (0,3) {$\bullet$};%
      \node (004) at (0,4) {$\bullet$};%
      \node (005) at (0,5) {$M_1$};
      \draw[>->] (000)--(001);%
      \draw[>->] (001)--(002);%
      \draw[>->] (002)--(003);%
      \draw[>->] (003)--(004);%
      \draw[>->] (004)--(005);%
      \begin{scope}[shift={(1,1)}]
        \node (100) at (0,0) {$\bullet$};%
        \node (101) at (0,1) {$\bullet$};%
        \node (102) at (0,2) {$\bullet$};%
        \node (103) at (0,3) {$\bullet$};%
        \node (104) at (0,4) {$N_1$};%
        \node (105) at (0,5) {$\bullet$};%
        \node (106) at (0,6) {$\bullet$};%
        \node (107) at (0,7) {$\bullet$};%
        \node (108) at (0,8) {$\bullet$};%
        \node (109) at (0,9) {$\bullet$};%
        \node (1010) at (0,10) {$M_2$};%
        \draw[>->] (100)--(101);%
        \draw[>->] (101)--(102);%
        \draw[>->] (102)--(103);%
        \draw[>->] (103)--(104);%
        \draw[>->] (104)--(105);%
        \draw[>->] (105)--(106);%
        \draw[>->] (106)--(107);%
        \draw[>->] (107)--(108);%
        \draw[>->] (108)--(109);%
        \draw[>->] (109)--(1010);%
      \end{scope}
      \begin{scope}[shift={(2,2)}]
        \node (200) at (0,0) {$\bullet$};%
        \node (201) at (0,1) {$\bullet$};%
        \node (202) at (0,2) {$\bullet$};%
        \node (203) at (0,3) {$\bullet$};%
        \node (204) at (0,4) {$\bullet$};%
        \node (205) at (0,5) {$\bullet$};%
        \node (206) at (0,6) {$\bullet$};%
        \node (207) at (0,7) {$\bullet$};%
        \node (208) at (0,8) {$\bullet$};%
        \node (209) at (0,9) {$\bullet$};%
        \draw[>->] (200)--(201);%
        \draw[>->] (201)--(202);%
        \draw[>->] (202)--(203);%
        \draw[>->] (203)--(204);%
        \draw[>->] (204)--(205);%
        \draw[>->] (205)--(206);%
        \draw[>->] (206)--(207);%
        \draw[>->] (207)--(208);%
        \draw[>->] (208)--(209);%
      \end{scope}
      \begin{scope}[shift={(3,3)}]
        \node (300) at (0,0) {$\bullet$};%
        \node (301) at (0,1) {$\bullet$};%
        \node (302) at (0,2) {$\bullet$};%
        \node (303) at (0,3) {$\bullet$};%
        \node (304) at (0,4) {$\bullet$};%
        \node (305) at (0,5) {$\bullet$};%
        \node (306) at (0,6) {$\bullet$};%
        \node (307) at (0,7) {$\bullet$};%
        \node (308) at (0,8) {$\bullet$};%
        \draw[>->] (300)--(301);%
        \draw[>->] (301)--(302);%
        \draw[>->] (302)--(303);%
        \draw[>->] (303)--(304);%
        \draw[>->] (304)--(305);%
        \draw[>->] (305)--(306);%
        \draw[>->] (306)--(307);%
        \draw[>->] (307)--(308);%
      \end{scope}
      \begin{scope}[shift={(4,4)}]
        \node (400) at (0,0) {$\bullet$};%
        \node (401) at (0,1) {$\bullet$};%
        \node (402) at (0,2) {$\bullet$};%
        \node (403) at (0,3) {$\bullet$};%
        \node (404) at (0,4) {$\bullet$};%
        \node (405) at (0,5) {$\bullet$};%
        \node (406) at (0,6) {$\bullet$};%
        \node (407) at (0,7) {$\bullet$};%
        \draw[>->] (400)--(401);%
        \draw[>->] (401)--(402);%
        \draw[>->] (402)--(403);%
        \draw[>->] (403)--(404);%
        \draw[>->] (404)--(405);%
        \draw[>->] (405)--(406);%
        \draw[>->] (406)--(407);%
      \end{scope}
      \begin{scope}[shift={(5,5)}]
        \node (500) at (0,0) {$S_c$};%
        \node (501) at (0,1) {$\bullet$};%
        \node (502) at (0,2) {$\bullet$};%
        \node (503) at (0,3) {$\bullet$};%
        \node (504) at (0,4) {$\bullet$};%
        \node (505) at (0,5) {$\bullet$};%
        \node (506) at (0,6) {$M_2'$};%
        \draw[>->] (500)--(501);%
        \draw[>->] (501)--(502);%
        \draw[>->] (502)--(503);%
        \draw[>->] (503)--(504);%
        \draw[>->] (504)--(505);%
        \draw[>->] (505)--(506);%
      \end{scope}
      \begin{scope}[shift={(6,6)}]
        \node (600) at (0,0) {$\bullet$};%
        \node (601) at (0,1) {$\bullet$};%
        \node (602) at (0,2) {$\bullet$};%
        \node (603) at (0,3) {$\bullet$};%
        \node (604) at (0,4) {$\bullet$};%
        \node (605) at (0,5) {$N_2$};%
        \node (606) at (0,6) {$\bullet$};%
        \node (607) at (0,7) {$\bullet$};%
        \node (608) at (0,8) {$\bullet$};%
        \node (609) at (0,9) {$\bullet$};%
        \node (6010) at (0,10) {$\cdots$};%
        \draw[>->] (600)--(601);%
        \draw[>->] (601)--(602);%
        \draw[>->] (602)--(603);%
        \draw[>->] (603)--(604);%
        \draw[>->] (604)--(605);%
        \draw[>->] (605)--(606);%
        \draw[>->] (606)--(607);%
        \draw[>->] (607)--(608);%
        \draw[>->] (608)--(609);%
        \draw[>->] (609)--(6010);%
      \end{scope}
      \begin{scope}[shift={(7,7)}]
        \node (700) at (0,0) {$\bullet$};%
        \node (701) at (0,1) {$\bullet$};%
        \node (702) at (0,2) {$\bullet$};%
        \node (703) at (0,3) {$\bullet$};%
        \node (704) at (0,4) {$\bullet$};%
        \node (705) at (0,5) {$\bullet$};%
        \node (706) at (0,6) {$\bullet$};%
        \node (707) at (0,7) {$\bullet$};%
        \node (708) at (0,8) {$\cdots$};%
        \node (709) at (0,9) {$\bullet$};%
        \draw[>->] (700)--(701);%
        \draw[>->] (701)--(702);%
        \draw[>->] (702)--(703);%
        \draw[>->] (703)--(704);%
        \draw[>->] (704)--(705);%
        \draw[>->] (705)--(706);%
        \draw[>->] (706)--(707);%
        \draw[>->] (707)--(708);%
        \draw[>->] (708)--(709);%
      \end{scope}
      \begin{scope}[shift={(8,8)}]
        \node (800) at (0,0) {$\bullet$};%
        \node (801) at (0,1) {$\bullet$};%
        \node (802) at (0,2) {$\bullet$};%
        \node (803) at (0,3) {$\bullet$};%
        \node (804) at (0,4) {$\bullet$};%
        \node (805) at (0,5) {$\bullet$};%
        \node (806) at (0,6) {$\cdots$};%
        \node (807) at (0,7) {$\bullet$};%
        \node (808) at (0,8) {$\bullet$};%
        \draw[>->] (800)--(801);%
        \draw[>->] (801)--(802);%
        \draw[>->] (802)--(803);%
        \draw[>->] (803)--(804);%
        \draw[>->] (804)--(805);%
        \draw[>->] (805)--(806);%
        \draw[>->] (806)--(807);%
        \draw[>->] (807)--(808);%
      \end{scope}
      \begin{scope}[shift={(9,9)}]
        \node (900) at (0,0) {$\bullet$};%
        \node (901) at (0,1) {$\bullet$};%
        \node (902) at (0,2) {$\bullet$};%
        \node (903) at (0,3) {$\bullet$};%
        \node (904) at (0,4) {$\cdots$};%
        \node (905) at (0,5) {$\bullet$};%
        \node (906) at (0,6) {$\bullet$};%
        \node (907) at (0,7) {$\bullet$};%
        \draw[>->] (900)--(901);%
        \draw[>->] (901)--(902);%
        \draw[>->] (902)--(903);%
        \draw[>->] (903)--(904);%
        \draw[>->] (904)--(905);%
        \draw[>->] (905)--(906);%
        \draw[>->] (906)--(907);%
      \end{scope}
      \begin{scope}[shift={(10,10)}]
        \node (1000) at (0,0) {$\bullet$};%
        \node (1001) at (0,1) {$\bullet$};%
        \node (1002) at (0,2) {$\cdots$};%
        \node (1003) at (0,3) {$\bullet$};%
        \node (1004) at (0,4) {$\bullet$};%
        \node (1005) at (0,5) {$\bullet$};%
        \node (1006) at (0,6) {$\bullet$};%
        \draw[>->] (1000)--(1001);%
        \draw[>->] (1001)--(1002);%
        \draw[>->] (1002)--(1003);%
        \draw[>->] (1003)--(1004);%
        \draw[>->] (1004)--(1005);%
        \draw[>->] (1005)--(1006);%
      \end{scope}
      \begin{scope}[shift={(11,11)}]
        \node (1100) at (0,0) {$\dots$};%
        \node (1101) at (0,1) {$\bullet$};%
        \node (1102) at (0,2) {$\bullet$};%
        \node (1103) at (0,3) {$\bullet$};%
        \node (1104) at (0,4) {$\bullet$};%
        \node (1105) at (0,5) {$\bullet$};%
        \draw[>->] (1100)--(1101);%
        \draw[>->] (1101)--(1102);%
        \draw[>->] (1102)--(1103);%
        \draw[>->] (1103)--(1104);%
        \draw[>->] (1104)--(1105);%
      \end{scope}
      \begin{scope}[shift={(12,12)}]
        \node (1200) at (0,0) {$\bullet$};%
        \node (1201) at (0,1) {$\bullet$};%
        \node (1202) at (0,2) {$\bullet$};%
        \node (1203) at (0,3) {$\bullet$};%
        \node (1204) at (0,4) {$N_d$};%
        \node (1205) at (0,5) {$M_{d+1}$};%
        \draw[>->] (1200)--(1201);%
        \draw[>->] (1201)--(1202);%
        \draw[>->] (1202)--(1203);%
        \draw[>->] (1203)--(1204);%
        \draw[>->] (1204)--(1205);%
      \end{scope}
      \begin{scope}[shift={(13,13)}]
        \node (1300) at (0,0) {$\bullet$};%
        \node (1301) at (0,1) {$\bullet$};%
        \node (1302) at (0,2) {$\bullet$};%
        \node (1303) at (0,3) {$\bullet$};%
        \node (1304) at (0,4) {$\bullet$};%
        \draw[>->] (1300)--(1301);%
        \draw[>->] (1301)--(1302);%
        \draw[>->] (1302)--(1303);%
        \draw[>->] (1303)--(1304);%
      \end{scope}
      \begin{scope}[shift={(14,14)}]
        \node (1400) at (0,0) {$\bullet$};%
        \node (1401) at (0,1) {$\bullet$};%
        \node (1402) at (0,2) {$\bullet$};%
        \node (1403) at (0,3) {$\bullet$};%
        \draw[>->] (1400)--(1401);%
        \draw[>->] (1401)--(1402);%
        \draw[>->] (1402)--(1403);%
      \end{scope}
      \begin{scope}[shift={(15,15)}]
        \node (1500) at (0,0) {$\bullet$};%
        \node (1501) at (0,1) {$\bullet$};%
        \node (1502) at (0,2) {$\bullet$};%
        \draw[>->] (1500)--(1501);%
        \draw[>->] (1501)--(1502);%
      \end{scope} 
      \begin{scope}[shift={(16,16)}]
        \node (1600) at (0,0) {$\bullet$};%
        \node (1601) at (0,1) {$\bullet$};%
        \draw[>->] (1600)--(1601);%
      \end{scope} 
      \begin{scope}[shift={(17,17)}]
        \node (1700) at (0,0) {$S_b$};%
      \end{scope}
      \draw[->>] (001)--(100);%
      \draw[->>] (002)--(101);%
      \draw[->>] (003)--(102);%
      \draw[->>] (004)--(103);%
      \draw[->>] (005)--(104);%
      
      \draw[->>] (101)--(200);%
      \draw[->>] (102)--(201);%
      \draw[->>] (103)--(202);%
      \draw[->>] (104)--(203);%
      \draw[->>] (105)--(204);%
      \draw[->>] (106)--(205);%
      \draw[->>] (107)--(206);%
      \draw[->>] (108)--(207);%
      \draw[->>] (109)--(208);%
      \draw[->>] (1010)--(209);%

      \draw[->>] (201)--(300);%
      \draw[->>] (202)--(301);%
      \draw[->>] (203)--(302);%
      \draw[->>] (204)--(303);%
      \draw[->>] (205)--(304);%
      \draw[->>] (206)--(305);%
      \draw[->>] (207)--(306);%
      \draw[->>] (208)--(307);%
      \draw[->>] (209)--(308);%

      \draw[->>] (301)--(400);%
      \draw[->>] (302)--(401);%
      \draw[->>] (303)--(402);%
      \draw[->>] (304)--(403);%
      \draw[->>] (305)--(404);%
      \draw[->>] (306)--(405);%
      \draw[->>] (307)--(406);%
      \draw[->>] (308)--(407);%

      \draw[->>] (401)--(500);%
      \draw[->>] (402)--(501);%
      \draw[->>] (403)--(502);%
      \draw[->>] (404)--(503);%
      \draw[->>] (405)--(504);%
      \draw[->>] (406)--(505);%
      \draw[->>] (407)--(506);%

      \draw[->>] (501)--(600);%
      \draw[->>] (502)--(601);%
      \draw[->>] (503)--(602);%
      \draw[->>] (504)--(603);%
      \draw[->>] (505)--(604);%
      \draw[->>] (506)--(605);%

      \draw[->>] (601)--(700);%
      \draw[->>] (602)--(701);%
      \draw[->>] (603)--(702);%
      \draw[->>] (604)--(703);%
      \draw[->>] (605)--(704);%
      \draw[->>] (606)--(705);%
      \draw[->>] (607)--(706);%
      \draw[->>] (608)--(707);%
      \draw[->>] (609)--(708);%
      \draw[->>] (6010)--(709);%

      \draw[->>] (701)--(800);%
      \draw[->>] (702)--(801);%
      \draw[->>] (703)--(802);%
      \draw[->>] (704)--(803);%
      \draw[->>] (705)--(804);%
      \draw[->>] (706)--(805);%
      \draw[->>] (707)--(806);%
      \draw[->>] (708)--(807);%
      \draw[->>] (709)--(808);%

      \draw[->>] (801)--(900);%
      \draw[->>] (802)--(901);%
      \draw[->>] (803)--(902);%
      \draw[->>] (804)--(903);%
      \draw[->>] (805)--(904);%
      \draw[->>] (806)--(905);%
      \draw[->>] (807)--(906);%
      \draw[->>] (808)--(907);%

      \draw[->>] (901)--(1000);%
      \draw[->>] (902)--(1001);%
      \draw[->>] (903)--(1002);%
      \draw[->>] (904)--(1003);%
      \draw[->>] (905)--(1004);%
      \draw[->>] (906)--(1005);%
      \draw[->>] (907)--(1006);%

      \draw[->>] (1001)--(1100);%
      \draw[->>] (1002)--(1101);%
      \draw[->>] (1003)--(1102);%
      \draw[->>] (1004)--(1103);%
      \draw[->>] (1005)--(1104);%
      \draw[->>] (1006)--(1105);%

      \draw[->>] (1101)--(1200);%
      \draw[->>] (1102)--(1201);%
      \draw[->>] (1103)--(1202);%
      \draw[->>] (1104)--(1203);%
      \draw[->>] (1105)--(1204);%

      \draw[->>] (1201)--(1300);%
      \draw[->>] (1202)--(1301);%
      \draw[->>] (1203)--(1302);%
      \draw[->>] (1204)--(1303);%
      \draw[->>] (1205)--(1304);%

      \draw[->>] (1301)--(1400);%
      \draw[->>] (1302)--(1401);%
      \draw[->>] (1303)--(1402);%
      \draw[->>] (1304)--(1403);%

      \draw[->>] (1401)--(1500);%
      \draw[->>] (1402)--(1501);%
      \draw[->>] (1403)--(1502);%

      \draw[->>] (1501)--(1600);%
      \draw[->>] (1502)--(1601);%

      \draw[->>] (1601)--(1700);%

      \draw [
      thick,
      decoration={
        brace,
        mirror,
        raise=0.5cm
      },
      decorate
      ] (0,0) -- (5,5) 
      node [pos=0.5,anchor=north,yshift=-0.55cm] {$\ell=\operatorname{len}(M_1)$};
    \end{tikzpicture}
  }%
\end{resizebox}
\end{center}
Consider now the Nakayama algebra $B=\kk Q_B/I$ with Gabriel quiver
\[
  \textstyle Q_B\colon\quad 1\to2\to\dots\to p,\qquad
  p=\sum_{i=1}^{d+1}\len{M_i}-\sum_{j=1}^d\len{N_j},
\]
and Kupisch series
\begin{align*}
  (&1,2,\dots,\len{M_1},\len{N_1}+1,\len{N_1}+2,\dots,\len{M_2},\\
   &\dots,\len{N_{d-1}}+1,\len{N_{d-1}}+2,\dots,\len{M_d},\len{M_{d+1}}).
\end{align*}
Observe that the above commutative diagram is indexed by the AR quiver of the
category $\mod{B}$ of finite-dimensional (right) $B$-modules. In other words, we
can extend the exact sequence \eqref{eq:delta} to an exact functor
$F\colon\mod{B}\to\AA$ that sends the simple $B$-modules to simple objects of
$\AA$. The exact functor $F$ lifts to an $\A$-functor
\[
  \widetilde{\F}\colon\dg{\Db{\mod{B}}}\longrightarrow\dg{\Db{\AA}}
\]
between the corresponding bounded derived DG categories, for example because
$\dg{\Db{\mod{B}}}$ is the pre-triangulated hull of $\mod{B}$~\cite[Thm.~6.1 and
Ex.~6.2]{Che23}. The many-objects version of the Homotopy Transfer Theorem
yields minimal models of $\dg{\Db{\mod B}}$ and $\dg{\Db{\AA}}$ so that we may
consider the induced $\A$-functor
\[
  \begin{tikzcd}
    H^*(\dg{\Db{\mod{B}}})\rar[dotted]{\F}\dar[swap]{\wr}&H^*(\dg{\Db{\AA}})\\
    \dg{\Db{\mod{B}}}\rar{\widetilde{\F}}&\dg{\Db{\AA}}\uar[swap]{\wr}
  \end{tikzcd}\qquad
  \begin{tikzcd}
    H^*(\dg{\Db{\mod{B}}})\rar[dotted]{\F_1}\dar[swap]{\mathrm{id}}&H^*(\dg{\Db{\AA}})\\
    H^*(\dg{\Db{\mod{B}}})\rar{H^*(\widetilde{\F})}&H^*(\dg{\Db{\AA}})\uar[swap]{\mathrm{id}}
  \end{tikzcd}
\]
We obtain, in particular, an $\A$-morphism
\[
  \F\colon\Ext[B]{S}{S}\longrightarrow\Ext{S}{S},
\]
where the Yoneda algebras are now endowed with minimal $\A$-algebra structures.

We claim that the underlying morphism of graded algebras
\begin{equation}
  \label{eq:F1(gamma)}
  \F_1\colon\Ext[B]{S}{S}\longrightarrow\Ext{S}{S},\qquad
  \gamma\longmapsto\delta,
\end{equation}
maps the Yoneda class $\gamma$ of an augmented minimal projective resolution of
the simple $B$-module concentrated at the vertex $p$ of the quiver $Q_B$ to the
class $\class{\delta}$. Indeed, $\gamma=\alpha_{d+1}\alpha_2\cdots\alpha_1$
is the product of Yoneda classes of short exact sequences between certain indecomposable
$B$-modules and, therefore,
\[
  \F_1(\gamma)=\F_1(\alpha_{d+1}\alpha_d\cdots\alpha_1)=\F_1(\alpha_{d+1})\F_1(\alpha_d)\cdots\F_1(\alpha_1)=\delta.
\]
Here, the classes $\F_1(\alpha_i)$, $i=1,\dots,d+1$, satisfy the equality on the
right-hand side by construction since the restriction of the graded functor
$\F_1\colon H^*(\dg{\Db{\mod{B}}})\to H^*(\dg{\Db{\AA}})$ to the graded
subcategories spanned by $\mod{B}$ and $\AA$, respectively, is induced by the
exact functor $F\colon\mod{B}\to\AA$ above.

Since Nakayama algebras are monomial algebras, an explicit description of a
minimal model of the Yoneda algebra $\Ext[B]{S}{S}$ is available due to
independent work of Chuang and King (unpublished) and of
Tamaroff~\cite[Thm.~4.9]{Tam21}.\footnote{The computation of the explicit
  product that we need is straightforward using the HTT; see also~\cite{Her21}.}
In particular, we can assume that\footnote{Compare with~\cite[Thm.~1]{Min16}
  where the equivalent of the Massey product
  $\langle\eta_1,\dots,\eta_{\ell-1},\eta\rangle\ni\pm\gamma$ appears.}
\begin{equation}
  \label{eq:m_ell}
  m_\ell^B(\eta_1\otimes\eta_2\otimes\dots\otimes\eta_{\ell-1}\otimes\eta)=\gamma,\qquad\ell=\len{M_1},
\end{equation}
where $\eta_i\in\Ext[B]<1>{S_{i+1}}{S_i}\cong\kk$ is the Yoneda class of an AR
sequence and $\eta\in\Ext[B]<d>{S_p}{S_\ell}\cong\kk$ is the Yoneda class of an
exact sequence of the form
\[
  0\to S_\ell\to M_2'\to\cdots\to M_{d+1}\to S_p\to0,
\]
with $M_2'$ the injective hull of the simple $B$-module $S_{\ell}$ and $M_{d+1}$
the projective cover of $S_p$, as indicated in the above diagram that we now
interpret as the AR quiver of $\mod{B}$. Moreover, an elementary computation
shows that
  \begin{align}
    \label{eq:m-vanishing}
    m_{j-i+1}^B(\eta_i\otimes \eta_{i+1}\otimes \cdots\otimes \eta_j)&\in\Ext[B]<2>{S_{j+1}}{S_i}=0
                                                                       \intertext{whenever $1\leq i<j<\ell$, and also}\label{eq:m-vanishing-x}
                                                                       m_{\ell-i+1}^B(\eta_i\otimes \eta_{i+1}\otimes \cdots\otimes \eta_{\ell-1}\otimes \eta)&\in\Ext[B]{S_p}{S_i}=0
  \end{align}
  whenever $1<i<\ell$. Here we use that the operation $m_k^B$ has degree $2-k$.

  Finally, the $\A$-morphism $\F\colon\Ext[B]{S}{S}\to\Ext{S}{S}$ satisfies in
  particular an equation of the form
  \[
    0=\partial(\F_\ell)=\sum_{\substack{r+1+t=k\\r+s+t=\ell}}\pm\F_k\circ(1^{\otimes
      r}\otimes m_s^{B}\otimes 1^{\otimes t})-\sum_{\substack{2\leq{k}\leq\ell\\
        i_1+\cdots+i_k=\ell}}\pm
    m_k^{\AA}\circ(\F_{i_1}\otimes\cdots\otimes\F_{i_k}),
  \]
  where $\partial=0$ since the $\A$-algebras involved are minimal. Therefore,
  using equations \eqref{eq:F1(gamma)} and \eqref{eq:m_ell},\footnote{Compare
    with~\cite[Thm.~A(ii)]{BMM20} and see also Section~3 therein.}
  \begin{align*}
    \delta=\F_1(\gamma)&=\F_1(m_\ell^B(\eta_1\otimes\cdots\otimes\eta_{\ell-1}\otimes\eta))
    \\&=m_\ell^{\AA}(\underbrace{\F_1(\eta_1)\otimes\cdots\otimes\F_1(\eta_{\ell-1})\otimes\F_1(\eta)}_{\in\Ext<1,d>{S}{S}})+\omega,
  \end{align*}
  where
  \[
    \omega=\sum_{\substack{r+1+t=k>1\\r+s+t=\ell}}\pm\F_k(\underbrace{\phantom{a}\cdots\phantom{a}}_{=0})+\sum_{\substack{2\leq{k}<\ell\\i_1+\cdots+i_k=\ell}}
    \pm
    m_k^{\AA}(\underbrace{\phantom{abc}\cdots\phantom{abc}}_{\in\Ext<1,d>{S}{S}}).
  \]
  To argue the vanishing of the first term we use that the inputs involve one
  the higher products \eqref{eq:m-vanishing}-\eqref{eq:m-vanishing-x} that we
  know vanish, and for the condition on the degrees of the inputs in the second
  term we note that
  \begin{align*}
    \F_{j-i+1}(\eta_i\otimes\cdots\otimes\eta_j)&\in\Ext<1>{S}{S}&&1\leq i<j<\ell,\\
    \F_{\ell-i+1}(\eta_i\otimes\cdots\otimes\eta_{\ell-1}\otimes\eta)&\in\Ext<d>{S}{S}&&1<i<\ell,
  \end{align*}
  since the component $\F_k$ is a morphism of degree $1-k$. This finishes the proof.
\end{proof}

\begin{remark}
  In a MathOverflow post~\cite{MadKel}, Madsen gives an outline of Keller's
  original proof of the theorem, which is
  unpublished.\footnote{But see~\cite[Thm.~3.22]{KM21a}.} Keller's proof uses the
  description of the category of finite-dimensional modules over a
  finite-dimensional algebra $A$ as the category of twisted stalks on the Yoneda
  algebra $\Ext[A]{S}{S}$ with its minimal $\A$-algebra structure. In fact, one
  may replace the ambient abelian category by the exact subcategory of objects
  filtered by a finite collection of objects, with $\operatorname{Ext}^*$ now
  understood in this exact subcategory. Our proof, which can be adapted to this
  more general setting with the appropriate modifications, avoids making explicit use of
  the category of twisted stalks by instead reducing the problem to a computation
  with Nakayama algebras.
\end{remark}

\subsection*{Acknowledgements}

The author thanks Julian K{\"u}lshammer for informing him of Madsen's
MathOverflow post~\cite{MadKel} and for interesting discussions on possible
extensions of Keller's theorem. The author also thanks Bernhard Keller for
clarifications on his original proof of the theorem.
\vspace{-0.5em}
\subsection*{Financial support}

The author's research is supported by the Swedish Research Council
(Vetenskapsrådet) Research Project Grant 2022-03748 `Higher structures in
higher-dimensional homological algebra.'

\printbibliography

\end{document}